\documentclass{article}

\usepackage{graphicx}
\usepackage{amsfonts}
\usepackage{amsmath}
\usepackage{amssymb}

\newtheorem{theorem}{Theorem}

\newtheorem{corollary}{Corollary}

\newtheorem{definition}{Definition}

\newtheorem{proposition}{Proposition}
\newtheorem{remark}{Remark}

\newenvironment{proof}[1][Proof]{\textbf{#1.} }{\ \rule{0.5em}{0.5em}}

\begin{document}

\title{Invariants of solvable rigid Lie algebras up to dimension $8$}

\author{Rutwig Campoamor-Stursberg\\Depto. Geometr\'{\i}a y Topolog\'{\i}a\\Fac. CC. Matem\'aticas U.C.M.\\
E-28040 Madrid ( Spain)\newline e-mail: rutwig@nfssrv.mat.ucm.es}

\date{}

\maketitle

\begin{abstract}
The invariants of all complex solvable rigid Lie algebras up to dimension eight are computed. Moreover we show, for rank one solvable algebras, some criteria to deduce to non-existence of non-trivial invariants or the existence of fundamental sets of invariants formed by rational functions of the Casimir invariants of the associated nilradical. 
\end{abstract}

\section{Introduction}

The importance of Casimir invariants of Lie algebras has increased its importance in representation theory and physics since the classical results due to Casimir and Racah \cite{Cas,Ra}. They allow for example a labelling of irreducible representations of Lie algebras, and in physical applications, the eigenvalues of invariant functions of the symmetry group  provide characterizations of a physical system, such as energy spectra. Casimir invariants of Lie algebras, which are elements of the centre of the universal enveloping algebra $\frak{A}$ of a Lie algebra, are commonly computed by using coadjoint orbits \cite{Di,AL}. For semi-simple Lie algebras the question has been solved long ago \cite{Ra}, but for solvable Lie algebras the question remains open. Here we obtain invariants which are not more Casimir invariants in the classical sense, since they may be transcendental functions. Their number and precise structure is also unknown, up to certain special classes which have been treated in recent literature \cite{NW2,NW,Nd,Pa}. 

In this paper we determine explicitely all invariant functions of complex solvable rigid Lie algebras up to dimension $8$. The interest of rigid structures is justified, for example, by the fact that the Lie algebra of the homogeneous Galilei group deforms onto the Lorentz algebra, which is rigid for being semi-simple \cite{Vi}. The classification of eight dimensional solvable rigid Lie algebras is a consequence of the root theory for solvable algebras developed in \cite{AG1}, and it allows theoretically the full determination of solvable rigid Lie algebras in any fixed dimension. Analyzing the structure of the toral subalgebras of such Lie algebras, we will be able to isolate properties that are valid in any dimension. For example, we will find sufficiency criteria to ensure the non-existence of nontrivial invariants for certain classes of rank one Lie algebras, or conditions to obtain the invariants of a solvable rigid Lie algebra as rational functions of the invariants of its nilradical. This points out that, for decomposable solvable Lie algebras, the key to the determination of invariants depends more on the way a torus acts on the nilradical than on the structure of the latter. This will also enable us to find explicitly the invariants of solvable rigid Lie algebras by extrapolation of the solutions in lower dimensional cases. This procedure is of interest, since it indicates the existence of certain rigid structures that depend heavily on their dimension. On the other side, nilradicals of rigid solvable Lie algebras are non-abelian, and their study enlarges the efforts made for these algebras in the last years \cite{NW2,Nd}. The rigid case allows us to show that for solvable Lie algebras the results will not be so harmonious as in the semi-simple case, since there are Lie algebras whose number of functionally independent invariants depends not on the dimension of a Cartan subalgebra, but on the dimension of the algebra. This is a succint indication that the non-decomposable case is, without doubt, the most difficult to solve. Probably the study of specific classes of Lie algebras is the only manner to approach the general solvable case, as classifications of these algebras only exist up to dimension six, and more advances seem not to be very probable ( among other reasons, because they imply the classification of nilpotent Lie algebras, classified only up to dimension seven). Resuming, the rigid case should be a first step to obtain interesting results on decomposable solvable Lie algebras, by studying the subtori of a maximal torus of derivations. 

\bigskip
Any Lie algebra considered in this work is finite dimensional over the field $\mathbb{C}$. Moreover, any $n$-dimensional Lie algebra $\frak{g}=\left(\mathbb{C}^{n},\mu\right)$ is identified with its law $\mu$ in the variety $\mathcal{L}^{n}$ \cite{AG3}. Moreover, we convene that nonwritten brackets are either zero or obtained by antisymmetry.

\section{Decomposability of rigid Lie algebras}

A Lie algebra $\frak{g}$ is a vector space $V$ endowed with an alternating bilinear product which satisfies the Jacobi identity. If $\mathcal{M}$ denotes the algebraic set consisting of all Lie algebra laws on $V$, then the linear group $GL\left(V\right)$ acts on the space of all alternating bilinear forms over $V$ by
\begin{eqnarray}
\left(f*\mu\right)=f^{-1}\left(\mu\left(f ,f \right)\right),\quad f\in GL(V), \mu\in \mathcal{O}(V)
\end{eqnarray}
showing that $\mathcal{M}$ is stable under this action. Therefore, the orbits $\mathcal{O}(\mu)$ of the linear group $GL\left(V\right)$ on $\mathcal{M}$ correspond to the isomorphism classes of Lie algebra laws on $V$. This allows to identify a Lie algebra $\frak{g}$ with the pair $\left(V,\mu\right)$, where $\mu$ is an element of the orbit $\mathcal{O}(\mu)$. We say that a Lie algebra $\frak{g}=\left(V,\mu\right)$ is rigid if the orbit $\mathcal{O}(\mu)$ is open in $\mathcal{M}$. Thus, roughly speaking, a Lie algebra $\frak{g}$ is rigid if any Lie algebra $\frak{g}^{\prime}$ close to $\frak{g}$ is isomorhic to $\frak{g}$. This led to the famous rigidity theorem of Nijenhuis and Richardson \cite{NR}, though, as pointed out by Richardson in \cite{Ri}, there are rigid Lie algebras whose second adjoint cohomology group is non-zero. As known, semi-simple Lie algebras have vanishing second adjoint cohomology groups, as a consequence of the Whitehead lemmas, and by the Levi decomposition theorem, the study of rigid Lie algebras reduces mainly to analyze solvable Lie algebras. Now rigidity imposes strong structural conditions in the solvable case, which simplify its description. The most important property in this sense is their decomposability.

\begin{definition}
A Lie algebra $\frak{g}$ is called decomposable if it can be written as
\[
\frak{g}=\frak{s}\oplus\frak{t}\oplus\frak{n}%
\]
where $\frak{s}$ is a Levi subalgebra, $\frak{n}$ the nilradical and
$\frak{t}$ an abelian subalgebra whose elements are $ad$-semisimple and which
satisfies $\left[  \frak{s+t},\frak{t}\right]  =0$.
\end{definition}

The abelian subalgebra $\frak{t}$ of $Der\frak{g}$ defined by
\[
\frak{t}=\{adX,\quad X\in\frak{t}\}
\]
is called, following Malcev \cite{Ma}, an exterior torus on $\frak{g}$. It is called
maximal torus, if it is maximal for the inclusion. As proven by Malcev, all maximal tori are pairwise conjugated, thus their dimension constitutes an invariant of the algebra, 
called the rank of $\frak{g}$ and noted by $r(\frak{g}).$

\bigskip
By the Levi decomposition theorem, we can restrict ourselves to analyze the solvable Lie algebras. Now semi-simple Lie algebras are rigid by the Nijenhuis-Richardson theorem \cite{NR}, so that it seems reasonable to begin with the study of solvable rigid Lie algebras in order to obtain criteria that are not dependent on the dimension. Now the structure of solvable rigid Lie algebras is known, and it simplifies the question is some manner. The central result is due to Carles \cite{Car}:

\begin{theorem}
A rigid Lie algebra $\frak{g}$ is algebraic. In particular, it is decomposable.
\end{theorem}

In particular, if the Levi subalgebra $\frak{s}$ is zero, then $\frak{g}$ is a solvable Lie algebra. In order to study such algebras, Ancochea and Goze developed a root theory \cite{AG1} for rigid solvable algebras. We recall this theory briefly: \newline 

Let $\frak{r}$ be a complex solvable decomposable Lie algebra, and 
$\frak{t}$ be a maximal torus.  Let $X$ be a non-zero vector such that
$ad_{\mu_{0}}X$ belongs to $\frak{t}$ .

\begin{definition}
\ We say that $X\in \frak{t}$ is regular if the dimension of 
\begin{equation*}
V_{0}\left( X\right) =\left\{ Y\in \frak{g}\quad |\quad [X,Y]=0\right\}
\end{equation*}
is minimal that is, $\dim V_{0}\left( X\right) \leq \dim V_{0}\left(
Z\right) $ for all $Z$ such that $adZ$ belongs to $\frak{t}$.
\end{definition}

Choose a regular vector $X$ and let be $p=\dim V_{0}\left( X\right) $. Consider a
basis ($X,Y_{1},...,Y_{n-p},$ $X_{1},...,X_{p-1})$ of eigenvectors of $adX$
such that $\left( X,X_{1},...,X_{p-1}\right) $ is a basis of $V_{0}\left(
X\right) $, $\left( Y_{1},...,Y_{n-p},X_{1},...,X_{k_{0}}\right) $ is a
basis of the maximal nilpotent ideal $\frak{n}$ of $\frak{g}$, and $\left(
X_{k_{0}+1},...,X_{p-1}\right) $ are vectors such that $adX_{i}\in $ $T$.

\begin{definition}
Suppose that $\frak{g}$ is not nilpotent. The root system of $\frak{g}$
associated to $\left( X,Y_{1},...,Y_{n-p},X_{1},...,X_{p-1}\right) $ is the
linear system (S) defined by the following equations :

$x_{i}+x_{j}=x_{k}$ if the $X_{k}$-component of $[X_{i},X_{j}]$ is
non-zero.

$y_{i}+y_{j}=y_{k}$ if the $Y_{k}$-component of $[Y_{i},Y_{j}]$ is
non-zero.

$x_{i}+y_{j}=y_{k}$ if the $Y_{k}$-component of $[X_{i},Y_{j}]$ is
non-zero.

$y_{i}+y_{j}=x_{k}$ if the $X_{k}$component of $[Y_{i},Y_{j}]$ is
non-zero.
\end{definition}

\begin{theorem}
If $rank\left( S\right) \neq \dim \left( \frak{n}\right) -1$, then $\frak{g}$ is not rigid.
\end{theorem}

See \cite{AG1} for a proof. There are two important corollaries which will be of importance for understanding the structure of solvable rigid laws.

\begin{corollary}
If $\frak{r}=\frak{n}\oplus\frak{t}$ is rigid, then $\frak{t}$ is a maximal torus over $\frak{n}$.
\end{corollary}

\begin{corollary}
If $\frak{g}$ is rigid then there is regular vector $X$ such that $ad_{\frak{g}}X$ is diagonal and its eigenvalues are integers.
\end{corollary}

These properties determine if a given Lie algebra is rigid or not. For
example, let us suppose that all elements of $V_{0}\left( X\right) $ are
semi-simple.\ If 
\begin{equation*}
rank\left( S\right) \neq \dim D^{1}\left( \frak{g}\right) -1
\end{equation*}
where $D^{1}\left( \frak{g}\right) $ is the derived subalgebra of $\frak{g}$%
, then $\frak{g}$ is not rigid.

\begin{remark}
Even if the roots can be choosen in $\mathbb{Z}$, this does not in general imply that the Lie algebra is rational. In \cite{AG3} various examples of this have been worked out.
\end{remark}

\section{Basic results on invariants. Notation}

As known, the dual representation $ad^{*}$ to the adjoint representation $ad$ of a Lie algebra $\frak{g}$ is called the coadjoint representation \cite{Di}. The problem of finding its invariants is indeed reduced to that of solving a system of linear first order partial differential equations. If $B=\left\{X_{1},..,X_{n}\right\}$ is a basis of the $n$-dimensional Lie algebra $\frak{g}$ and $\left\{x_{1},..,x_{n}\right\}$ a coordinate system on the dual space, then the infinitesimal generators of the action are denoted by $\widetilde{X}_{i}$. If moreover the structure constants of $\frak{g}$ are given by $[X_{i},X_{j}]=C_{i,j}^{k}X_{k}$ over the basis $B$, the it follows that a function $F\in C^{\infty}\left(\frak{g}^{*}\right)$ is an invariant of the coadjoint representation if and only if the two following conditions are satisfied:
\begin{enumerate}

\item $\widetilde{X}_{i}=\sum_{j}(-C_{i,j}^{k})x_{k}\frac{\partial}{\partial{x_{j}}}$ and $\left[\widetilde{X}_{i},\widetilde{X}_{j}\right]=C_{i,j}^{k}\widetilde{X}_{k}$,

\item $F$ satisfies the system of linear first order PDE $\widetilde{X}_{i}F=0,\quad 1\leq i\leq n$
\end{enumerate}

Solutions to this system are usually found by integration of the corresponding system of characteristic equations or other standard integration procedures \cite{D}.A maximal set of functionally independent solutions will be called a fundamental set of invariants. The determination of the invariants as solutions of the associated system of linear first order PDEs is the most common method to determine the invariants of Lie algebras \cite{AL,Pa,Wi,An}, and the one we will apply here. The cardinal of a fundamental set of invariants is $N=dim\left(\frak{g}\right)-r$, where $r$ is the rank of the commutator table considered as a matrix \cite{Be}. Since a Lie algebra law is an alternated tensor of type $(2,1)$, this rank does not depend on the basis chosen. Moreover, by antisymmetry, this rank is even, and we obtain that $N\equiv dim\left(\frak{g}\right) ( mod 2 )$. In particular, an odd dimensional Lie algebras has nontrivial invariants. From a result due to Dixmier \cite{Di,AL}, it is known that if the Lie algebra is algebraic, then we can find a maximal set of functionally independent solutions formed by rational functions. In view of Carles' algebraicity theorem, we in particular obtain the following result.

\begin{proposition}
A rigid Lie algebra $\frak{n}\oplus\frak{t}$ having non-trivial invariants admits a fundamental set of invariants formed by rational functions. 
\end{proposition}

The converse is easily seen to be false. For example, take the non-nilpotent Lie algebra $N_{6,1}^{\alpha \beta \gamma \delta}$ of \cite{Nd} with brackets $[X_{1},N_{1}]=\alpha N_{1}, [X_{1},N_{2}]=\gamma N_{2}, [X_{1},N_{4}]=N_{4}, [X_{2},N_{1}]=\beta N_{1}, [X_{2},N_{2}]=\delta N_{2}, [X_{2},N_{3}]=N_{3}$. It is easy to verify that it is non-rigid, as the torus generated by ${X_{1},X_{2}}$ is non-maximal, since the nilradical of  $N_{6,1}^{\alpha \beta \gamma \delta}$ is abelian of rank four, and has a fundamental set of invariants formed by the rational functions $\left\{\frac{n_{3}^{\beta}n_{4}^{\alpha}}{n_{1}}, \frac{n_{3}^{(\beta \gamma - \alpha \delta)}}{n_{1}^{\gamma}}\right\}$.

\section{Some sufficiency criteria}
In this paragraph we analyze some results which allow to deduce the inexistence of nontrivial invariants of a solvable rigid law $\frak{r}=\frak{n}\oplus\frak{t}$ by analyzing the structure of the torus $\frak{t}$.

\begin{proposition}
Let $dim\left(\frak{n}\oplus\frak{t}\right)=2n$. Suppose that the rank of $\frak{n}$ is one and that the nilradical $\frak{n}$ has a one dimensional center. If there exists a basis $\left\{Y_{1},..,Y_{2n}\right\}$ such that 
\begin{enumerate}
\item $\frak{t}$ is generated by $Y_{n}$,

\item $Y_{2n}\in Z\left(\frak{n}\right)$, and

\item $[Y_{j},Y_{2n-j}]=Y_{2n},\quad 2\leq j\leq n-1$
\end{enumerate}
then $\frak{r}=\frak{n}\oplus\frak{t}$ has no nontrivial invariants.
\end{proposition}

\begin{corollary}
For $n\geq 3$, any solvable Lie algebra $\frak{r}=\frak{n}\oplus\frak{t}$ whose torus $\frak{t}$ has the eigenvalues $\left(1,2,..,n-1,n+1,..,2n\right)$
has only trivial invariant functions.
\end{corollary}

In \cite{Nd} the author points out the importance of finding characterizations of those solvable Lie algebras with non-abelian nilradical  non admitting non-trivial invariant functions. Even in the decomposable case such a characterization seems not to be realizable, as the following example shows: Consider the nilpotent Lie algebra $\frak{n}$ of dimension seven given by the brackets $[X_{1},X_{2}]=X_{4}, [X_{1},X_{3}]=X_{5}, [X_{1},X_{6}]=X_{7}, [X_{2},X_{3}]=X_{6}, [X_{2},X_{5}]=X_{7}, [X_{3},X_{4}]=-X_{7}$. This algebra is of rank three. Take the one dimensional tori $\frak{t}_{1}$ and $\frak{t}_{2}$ whose action over the basis $\left\{X_{1},..,X_{7}\right\}$ is, respectively, given by the sequences $\left(1,0,0,1,1,0,1\right)$ and $\left(1,-1,0,0,1,-1,0\right)$. Both solvable non-nilpotent Lie algebras $\frak{n+t_{1}}$ and $\frak{n+t_{2}}$ have the same non-abelian nilradical, but the first algebra has no non-trivial invariants, while the second admits the polynomial invariant $y_{7}$, since its center is nonzero. This shows that, even if characterizations of solvable Lie algebras having only trivial solutions exist, they cannot be based on the structure of the nilradical, but on the way the torus (in the decomposable case) acts on this nilradical.\newline This example also provides an obvious criterion to guarantee the existence of polynomial solutions of a solvable Lie algebra:

\begin{proposition}
Let $\frak{n}$ be a nilpotent Lie algebra of rank $r\geq 1$ and $\frak{t}$ a maximal torus over $\frak{n}$. If there exists a toral subalgebra $\frak{t}^{\prime}\subset\frak{t}$ such that the action of $\frak{t}^{\prime}$ on a central ideal $I$ is zero, then the solvable Lie algebra $\frak{n}\oplus\frak{t}^{\prime}$ has at least $dim\left(I\right)$ polynomial solutions. Moreover, if the ideal $I$ is generated by $\left\{Y_{1},..,Y_{t}\right\}$, these solution can be chosen as  $\left\{y_{1},..,y_{t}\right\}$.
\end{proposition}

\section{Solvable rigid Lie algebras up to dimension six}

The determination of solvable rigid laws in dimensions up to six does not provide much information about the variety, since the number of these laws is very reduced. In dimension $2$ and $4$ ( there are no solvable rigid laws in dimension $3$ \cite{AG3}) the solvable rigid algebras are simply the semidirect product of abelian Lie algebras with a maximal torus, and, by the results of \cite{NW2, Wi}, have no nontrivial invariants. For dimension $5$ we have only the law $\frak{r}=\frak{h}_{1}\oplus \frak{t}$, where $\frak{h}_{1}$ is the $3$-dimensional Heisenberg Lie algebra. Here we find a rational invariant, which also follows from the analysis undertaken in \cite{Wi}. Finally, there are three solvable rigid laws in dimension $6$. As known, invariants of solvable Lie algebras up to this dimension have been determined in \cite{Nd,Pa}.

\begin{table}
\caption{\label{rigid6}Invariants of solvable rigid Lie algebras in dimension $6$.}
\begin{tabular}{@{}lll}
\bf
Algebra&{\rm Brackets}&{\rm Invariants}\\\hline
\rm
$\frak{r}_{6}^{1}$& $[V_{1},Y_{i}]=iY_{i},\quad i=1,2,3,4,5$ & none\\
& $[Y_{1},Y_{i}]=Y_{i+1},\quad i=2,3,4$ &\\
& $[Y_{2},Y_{3}]=Y_{5}$ &\\\hline
$\frak{r}_{6}^{2}$ & $[V_{1},Y_{i}]=iY_{i},\quad i=1,2,3,4$& none\\
& $[V_{2},Y_{i}]=Y_{i},\quad i=2,3,4$& \\
& $[Y_{1},Y_{i}]=Y_{i+1},\quad i=2,3$ &\\\hline
$\frak{r}_{6}^{3}$ & $[V_{1},Y_{i}]=iY_{i},\quad i=1,2,3$& none\\
& $[V_{2},Y_{2}]=Y_{2},$& \\
& $[V_{3},Y_{3}]=Y_{3}$ &\\\hline
\bf
\end{tabular}
\end{table}

\section{Solvable rigid Lie algebras in dimension $7$}
As follows from the parity of the dimension, seven dimensional rigid solvable Lie algebras have at least one nontrivial invariant function. It can easily be seen that any such algebra has at most rank three, and that their nilradical is non-abelian. Their classification follows easily from the classification of low dimensional Lie algebras \cite{Se}. \newline In keeping the notation used in previous sections, we illustrate the determination of the invariants of a solvable Lie algebra considering the rank one Lie algebra $\frak{r}_{7}^{3}$ given by the brackets $[V_{1},Y_{i}]=iY_{i}, i=1,3,4,5,6,7;\quad [Y_{1},Y_{i}]=Y_{i+1}, i=3,4,5,6;\quad [Y_{3},Y_{4}]=Y_{7}$ over the basis $\left\{Y_{1},Y_{3},Y_{4},Y_{5},Y_{6},Y_{7},V_{1}\right\}$. It is trivial to see that the system $\left\{\widetilde{Y}_{i}.F=0,i=1,3,4,5,6,7;\quad \widetilde{V}_{1}.F=0\right\}$ reduces to the two following equations:
\begin{eqnarray}
\widetilde{X_{1}^{\prime}}.F=(y_{6}\partial_{y_{5}}+y_{7}\partial_{y_{6}}).F=0\\
\widetilde{V_{1}^{\prime}}.F=(5y_{5}\partial_{y_{5}}+6y_{6}\partial_{y_{6}}+7y_{7}\partial_{y_{7}}).F=0
\end{eqnarray}
From the characteristic equations of $(6.1)$ we easily obtain the set of invariants $\left\{y_{7},y_{6}^{2}-2y_{5}y_{7}\right\}$. Since the reduced system is complete, it suffices to apply the usual procedures of systems of two equations in three variables \cite{D}: taking $u=y_{7}$ and $v= y_{6}^{2}-2y_{5}y_{7}$ we obtain $\widetilde{V_{1}^{\prime}}\left(u\right)=7u$, $\widetilde{V_{1}^{\prime}}\left(v\right)=12v$ and the equation
\begin{eqnarray}
\frac{\partial{F(u,v)}}{\partial{u}}+\frac{12v}{7u}\frac{\partial{F(u,v)}}{\partial{v}}=0
\end{eqnarray}
which provides the solution $\left\{\frac{v^{7}}{u^{12}}=\frac{(2y_{5}y_{7}-y_{6}^{2})^{7}}{y_{7}^{12}}\right\}$. As the rank of the commutator table ( interpreted as matrix ) is six, the preceding solution constitutes a fundamental set of invariants for $\frak{r}_{7}^{3}$.

\begin{table}
\caption{\label{rigid7}Invariants of solvable rigid Lie algebras in dimension $7$.}
\begin{tabular}{@{}lll}
\bf
Algebra&{\rm Brackets}&{\rm Invariants}\\\hline
\rm
$\frak{r}_{7}^{1}$ & $[V_{1},Y_{i}]=iY_{i},\quad i=1,2,3,4,5,6$& $I_{1}=\frac{(3y_{6}^{2}y_{3}+y_{5}^{3}-3y_{4}y_{5}y_{6})^{2}}{y_{6}^{5}}$ \\
 & $[Y_{1},Y_{i}]=Y_{i+1},\quad i=2,3,4,5$  &\\
 & $[Y_{2},Y_{i}]=Y_{i+2},\quad i=3,4.$ & \\\hline
$\frak{r}_{7}^{2}$ & $[V_{1},Y_{i}]=iY_{i},\quad i=1,2,3,4,5,7 $ & $I_{1}=\frac{(6y_{3}y_{5}y_{7}-2y_{5}^{3}-3y_{7}y_{4}^{2}+6y_{7}^{2}y1)^{7}}{y_{7}^{15}}$ \\
 & $[Y_{1},Y_{i}]=Y_{i+1,}\quad i=2,3,4$ & \\
 & $[Y_{2},Y_{i}]=Y_{i+2},\quad i=3,5 $& \\
 & $[Y_{3},Y_{4}]=-Y_{7}$ &\\\hline
$\frak{r}_{7}^{3}$& $[V_{1},Y_{i}]=iY_{i},\quad i=1,3,4,5,6,7$ & $I_{1}=\frac{(2y_{5}y_{7}-y_{6}^{2})^{7}}{y_{7}^{12}}$\\
 & $[Y_{1},Y_{i}]=Y_{i+1},\quad i=3,4,5,6,$ &\\
 & $[Y_{3},Y_{4}]=Y_{7}$ & \\\hline
$\frak{r}_{7}^{4}$ & $[V_{1},Y_{i}]=iY_{i},\quad i=1,2,3,4,5$ & $I_{1}=\frac{(2y_{5}y_{3}^{\prime}+y_{4}^{2}-2y_{3}y_{5})^{5}}{y_{5}^{8}}$\\
 & $[V_{1},Y_{3}^{\prime}]=3Y_{3}^{\prime}$, & \\
 & $[Y_{1},Y_{i}]=Y_{i+1},\quad i=2,3,4,$ & \\
 & $[Y_{2},Y_{3}]=[Y_{2},Y_{3}^{\prime}]=Y_{5}.$ & \\\hline
$\frak{r}_{7}^{5}$ & $[V_{1},Y_{i}]=iY_{i}, i=1,3,\quad [V_{1},V_{i}]=(i-2)Y_{i}, i=4,5$ & $I_{1}=\frac{(3y_{2}y_{5}^{2}+y_{4}^{3}-3y_{3}y_{4}y_{5})^{2}}{(2y_{3}y_{5}-y_{4}^{2})^{3}}$\\
 & $[ V_{2},Y_{i}]=Y_{i},\quad i=2,3,4,5,$ & \\
 & $[ Y_{1},Y_{i}]  =Y_{i+1},\quad i=2,3,4.$ & \\\hline
$\frak{r}_{7}^{6}$ & $[ V_{1},Y_{i}]  =iY_{i},\quad i=1,2,3,4,5,$ & $I_{1}=\frac{y_{4}^{2}(2y_{1}y_{5}+y_{3}^{2}-2y_{2}y_{4})^{3}}{y_{5}^{2}}  $\\
 & $[ V_{2},Y_{i}]=Y_{i},\quad i=2,3,4,$ &\\
 & $[ V_{2},Y_{5}]=2Y_{5},$ & \\
 & $[ Y_{1},Y_{i}]=Y_{i+1},\quad i=2,3,$ & \\
 & $[ Y_{2},Y_{3}]  =Y_{5}.$ & \\\hline
$\frak{r}_{7}^{7}$ & $[ V_{1},Y_{i}]=Y_{i},i=1,3,\quad [V_{1},V_{i}]=2Y_{i}, i=2,5$ & $I_{1}=\frac{-2y_{5}^{2}v_{2}+v_{1}y_{5}^{2}+2y_{2}y_{3}y_{5}+y_{1}y_{4}y_{5}-y_{2}y_{4}^{2}}{y_{5}^{2}}$\\
 & $[ V_{2},Y_{i}]=Y_{i},\quad i=3,4,5,$ & \\
 & $[ Y_{1},Y_{i}]=Y_{i+1},\quad i=3,4,$ & \\
 & $[ Y_{2},Y_{3}]  =Y_{5}$  & \\\hline
$\frak{r}_{7}^{8}$ & $[ V_{1},Y_{i}]  =iY_{i},i=1,3$ & $I_{1}=\frac{v_{1}y_{3}-v_{2}y_{3}+y_{1}y_{2}}{y_{3}}$\\
 & $[ V_{2},Y_{i}]  =Y_{i},\quad i=2,3,$ & \\
 & $[ V_{2}^{\prime},Y_{4}]=Y_{4},$ & \\
 & $[ Y_{1},Y_{2}]  =Y_{3}.$ & \\\hline
\bf
\end{tabular}
\end{table}

\begin{proposition}
Any solvable rigid Lie algebra $\frak{r}$ of dimension $7$ and rank one has a fundamental set of invariants formed by a quotients of invariant polynomials of the nilradical.
\end{proposition}

The proof is straightforward, as follows from table 2. Now this result can be generalized to arbitrary dimension, under an additional assumption. From the analysis of the existing lists of rigid Lie algebras, it has been conjectured in the middle 80's that any solvable rigid Lie algebra $\frak{r}$ has necessarily a trivial center \cite{Car,AG2,AC,Car2}. For rank $r\geq 2$ this fact is easily proven by deforming the torus, while for rank one the question remains a conjecture. The difficulty of this case is deeply related with the classification of complete Lie algebras and cohomological problems \cite{Car}.\newline In any case, if $\frak{r}=\frak{n}\oplus\frak{t}$ is a solvable rigid Lie algebra of rank one and   $X$ a central element of the nilradical $\frak{n}$, the infinitesimal generator $\widetilde{X}$ is of the form $\widetilde{X}=[T,X]\frac{\partial}{\partial{t}}$, and therefore an invariant function $F$ does not depend on $t$. The following proposition is an easy provable consequence of the structure of the remaining system $\widetilde{X}_{i}F=0$, since this system is formed by the linear system of PDE which gives the invariants of the nilradical, to which the equation which describes the action of the torus over $\frak{n}$ must be added: 

\begin{proposition}
Let $\frak{r}=\frak{n}\oplus \frak{t}$ be a solvable rigid Lie algebra of rank one and trivial center. Then $\frak{r}$ admits a fundamental set of invariants formed by quotients of invariant polynomials of the nilradical $\frak{n}$.
\end{proposition}

\section{Solvable rigid Lie algebras in dimension $8$}

The first classification of eight dimensional solvable rigid Lie algebras appeared in 1986, though the corrected complete list has recently been published in \cite{AG2}. It is based on the study of the eigenvalues of a regular operator $ad(X)$, and is also related with the classification of seven dimensional nilpotent Lie algebras \cite{Se}. The rank theorem allowed to complete the list and to detect the absence of two laws. The list we present here is basically the one given in \cite{AG2}, up to minor changes relative to the choice of regular vectors. This is also the first dimension where even dimensional solvable rigid laws admit nontrivial invariants. It is convenient to separate these algebras by rank.

\begin{table}
\caption{\label{rigid81}Invariants of solvable rigid Lie algebras in dimension $8$ and rank $1$.}
\begin{tabular}{@{}lll}
\bf
Algebra&{\rm Brackets}&{\rm Invariants}\\\hline
\rm
$\frak{r}_{8}^{1}$ & $[ V_{1},Y_{i}]  =iY_{i},\quad i=1,2,3,5,6,7,8,$ & none\\
 & $[ Y_{1},Y_{i}]=Y_{i+1},\ i=2,6,7$ & \\
 & $[ Y_{2},Y_{3}]=Y_{5}, \quad [ Y_{2},Y_{5}]=Y_{7},$ &\\
 & $[ Y_{2},Y_{6}]=Y_{8}, \quad [Y_{3},Y_{5}]=Y_{8}.$ &\\\hline
$\frak{r}_{8}^{2}$ & $[ V_{1},Y_{i}]=iY_{i},\quad i=1,3,4,5,6,7,8$ & $I_{1}=\frac{y_{8}^{7}}{(y_{7}^{2}-2y_{6}y_{8})^{4}}$\\ 
 & $[ Y_{1},Y_{i}]=Y_{i+1},\quad i=3,4,5,6,7$ & $I_{2}=\frac{(y_{6}^{2}-2y_{5}y_{7}+2y_{4}y_{8})^{2}}{y_{8}^{3}}$\\
 & $[ Y_{3},Y_{i}]=Y_{i+3},\quad i=4,5.$ &\\\hline
$\frak{r}_{8}^{3}$ & $[ V_{1},Y_{i}]  =iY_{i},\quad i=1,3,4,5,6,7,9$ & $I_{1}=\frac{y_{9}^{7}}{y_{7}^{9}}$\\
 & $[ Y_{1},Y_{i}]=Y_{i+1},\quad i=3,4,5,6$ & $I_{2}=\frac{(2y_{9}^{2}y_{1}-y_{6}^{2}y_{7}+2y_{4}y_{6}y_{9}-y_{5}^{2}y_{9}+2y_{7}^{2}y_{5}-2y_{3}y_{7}y_{9})^{7}}{y_{7}y_{9}^{7}}$\\
 & $[ Y_{3},Y_{i}]=Y_{i+3},\quad i=4,6,$ &\\
 & $[ Y_{4},Y_{5}]=-Y_{9}.$ &\\\hline
$\frak{r}_{8}^{4}$ & $[ V_{1},Y_{i}]=iY_{i},\quad i=1,4,5,6,7,8,9$ & $I_{1}=\frac{(y_{8}^{2}-2y_{7}y_{9})^{9}}{y_{9}^{16}}$\\
 & $[ Y_{1},Y_{i}]=Y_{i+1},\quad i=4,5,6,7,8$ & $I_{2}=\frac{y_{9}^{8}}{(3y_{6}y_{9}^{2}+y_{8}^{3}-3y_{7}y_{8}y_{9})^{3}}$\\
 & $[ Y_{4},Y_{5}]=Y_{9}.$ & \\\hline
$\frak{r}_{8}^{5}$ & $[ V_{1},Y_{i}]  =iY_{i},\quad i=2,3,4,5,6,7,8$ & $I_{1}=\frac{y_{7}^{8}}{y_{7}^{8}}$ \\
 & $[ Y_{2},Y_{i}]=Y_{i+2},\quad i=3,4,5,6 $ & $I_{2}=\frac{(2y_{8}^{2}y_{4}-y_{6}^{2}y_{8}+2y_{6}y_{7}^{2}-2y_{5}y_{7}y_{8})^{7}}{y_{7}^{4}y_{8}^{14}}$\\
 & $[ Y_{3},Y_{i}]  =Y_{i+3},\quad i=4,5.$ & \\\hline
$\frak{r}_{8}^{6}$ & $[ V_{1},Y_{i}]=iY_{i},\quad i=2,3,4,6,7,8,10$ & none \\
 & $[ Y_{2},Y_{i}]=Y_{2+i},\quad i=4,6,8,$ & \\
 & $[ Y_{3},Y_{i}]  =Y_{i+3},\quad i=4,7,$ & \\
 & $[ Y_{4},Y_{6}]  =Y_{10}.$ & \\\hline
$\frak{r}_{8}^{7}$ & $[ V_{1},Y_{i}]  =iY_{i},\quad i=2,3,5,6,7,8,9$ & $I_{1}=\frac{y_{9}^{8}}{y_{8}^{9}}$\\
 & $[ Y_{2},Y_{i}]=Y_{i+2},\quad i=3,5,6,7,$ & \\
 & $[ Y_{3},Y_{i}] =Y_{i+3},\quad i=5,6.$ & $I_{2}=\frac{(2y_{5}y_{9}^{2}-y_{7}^{2}y_{9}+2y_{7}y_{8}^{2}-2y_{6}y_{8}y_{9})^{8}}{y_{8}^{5}y_{9}^{16}}$ \\\hline
$\frak{r}_{8}^{8}$ &$[ V_{1},Y_{i}]  =iY_{i},\quad i=2,3,5,7,8,9,11$ & $I_{1}=\frac{y_{11}^{18}}{(y_{9}^{2}-2y_{7}y_{11})^{11}}$\\
 & $[ Y_{2},Y_{i}]=Y_{i+2},\quad i=3,5,7,9$ & $I_{2}=\frac{(y_{9}^{2}-2y_{7}y_{11})^{3}}{(3y_{5}y_{11}^{2}+y_{9}^{3}-3y_{7}y_{9}y_{11}-\frac{3}{2}y_{8}^{2}y_{11})^{2}}$\\
 & $[ Y_{3},Y_{i}]=Y_{i+3},\quad i=5,8.$ & \\\hline
$\frak{r}_{8}^{9}$ & $[ V_{1},Y_{i}]=iY_{i},\quad i=1,2,3,4,5,6$ & none\\
 &$[ V_{1},Y_{3}^{\prime}]=3Y_{3}^{\prime}$ & \\
 & $[ Y_{1},Y_{i}]  =Y_{i+1},\quad i=2,3,4,5,$ & \\
& $[ Y_{2},Y_{i}]  =Y_{i+2},\quad i=3,4,$ & \\
& $[ Y_{2},Y_{3}^{\prime}]=Y_{5},\quad [Y_{3},Y_{3}^{\prime}]=Y_{6}$ & \\\hline
$\frak{r}_{8}^{10}$& $[ V_{1},Y_{i}]  =iY_{i},\quad i=1,2,3,4,5,6,$ & $I_{1}=\frac{(2y_{6}y_{4}^{\prime}+y_{5}^{2}-2y_{4}y_{6})^{3}}{y_{6}^{5}}$\\
& $[ V_{1},Y_{4}^{\prime}]  =4Y_{4}^{\prime}$ & $I_{2}=\frac{(3y_{6}^{2}y_{3}+y_{5}^{3}-3y_{4}y_{5}y_{6})^{2}}{y_{6}^{5}}$\\
& $[ Y_{1},Y_{i}]  =Y_{i+1},\quad i=2,3,4,5$ & \\
& $[ Y_{2},Y_{i}]  =Y_{i+2},\quad i=3,4,$ &\\
& $[ Y_{2},Y_{4}^{\prime}]  =Y_{6}$ &\\\hline
\bf
\end{tabular}
\end{table}

As follows from table 3., any solvable rigid algebra of rank one, up to the algebras $\frak{r}_{8}^{1},\frak{r}_{8}^{6}$ and $\frak{r}_{8}^{9}$, has two invariants. For the first two ones this follows at once from proposition 5., while for the third the reasoning is similar.

\begin{table}
\caption{\label{rigid81} continued}
\begin{tabular}{@{}lll}
\bf
Algebra&{\rm Brackets}&{\rm Invariants}\\\hline
\rm
$\frak{r}_{8}^{11}$& $[ V_{1},Y_{i}]  =iY_{i},\quad i=1,2,3,4,5,6,$ & $I_{1}=\frac{(y^{\prime}_{5}-y_{5})^{6}}{y_{6}^{5}}$\\
& $[ V_{1},Y_{5}^{\prime}]=5Y_{5}^{\prime},\quad [ Y_{1},Y_{5}^{\prime}]  =Y_{6}$ & $I_{2}=\frac{(6y_{6}^{2}y_{3}+2y_{5}^{3}+3y_{5}^{2}(y_{5}^{\prime}-y_{5})-6y_{4}y_{5}y_{6}-6y_{4}(y_{5}^{\prime}-y_{5})y_{6})^{2}}{y_{6}^{5}}$\\
& $[ Y_{1},Y_{i}]  =Y_{i+1},\quad i=2,3,4,5,$ & \\
& $[ Y_{2},Y_{3}]=Y_{5}^{\prime},\quad [Y_{2},Y_{4}]=Y_{6}$ & \\\hline
$\frak{r}_{8}^{12}$& $[ V_{1},Y_{i}]  =iY_{i},\quad i=1,2,3,4,5,7,$ & $I_{1}=\frac{(2y_{3}y_{7}-y_{5}^{2})^{7}}{y_{7}^{16}}$\\
& $[ V_{1},Y_{3}^{\prime}]  =3Y_{3}^{\prime}, [ Y_{1},Y_{2}]  =Y_{3}^{^{\prime}}$ & $I_{2}=\frac{y_{7}^{30}(y_{5}^{4}y_{7}^{2}+2y_{3}y_{5}^{2}y_{7}^{3}+4y_{3}^{2}y_{7}^{4})^{7}}{(y_{7}^{2}y_{5}^{10}-4y_{3}y_{7}^{3}y_{5}^{8}+4y_{3}^{2}y_{5}^{6}y_{7}^{4}-8y_{3}^{3}y_{5}^{4}y_{7}^{5}+32y_{3}^{4}y_{5}^{2}y_{7}^{6}-32y_{3}^{5}y_{7}^{7})^{7}}$\\
& $[Y_{1},Y_{3}^{\prime}]=Y_{4}, [ Y_{1},Y_{4}]  =Y_{5}$ & \\
& $[Y_{2},Y_{3}]=Y_{5}, [Y_{2}^{^{\prime}},Y_{3}^{^{\prime}}]=Y_{5}$ & \\
& $[Y_{2},Y_{5}]=[ Y_{3}^{^{\prime}},Y_{4}]=Y_{7}.$ & \\\hline
$\frak{r}_{8}^{13}$ & $[ V_{1},Y_{i}]=iY_{i},\quad i=1,2,3,4,5,7,$ & $I_{1}=\frac{(y_{5}-y_{5}^{\prime})^{7}}{y_{7}^{5}}$\\
& $[ V_{1},Y_{5}^{\prime}]  =5Y_{5}^{\prime},$ & $I_{2}=\frac{(6y_{7}^{2}y_{1}-2y_{5}^{3}+3y_{5}^{2}(y_{5}-y_{5}^{\prime})+6y_{3}y_{5}y_{7}-3y_{7}y_{4}^{2})^{7}}{y_{7}^{15}}$ \\
& $[ Y_{1},Y_{i}]  =Y_{i+1},\quad i=2,3,4,$ & \\
&$[ Y_{2},Y_{3}]=Y_{5}^{\prime}, \quad [ Y_{3},Y_{4}]
=-Y_{7}$ & \\
&$[ Y_{2},Y_{5}]=[ Y_{2},Y_{5}^{\prime}]=Y_{7}$, & \\\hline
$\frak{r}_{8}^{14}$ & $[ V_{1},Y_{i}]=iY_{i},\quad i=1,3,4,5,6,7,$ & $I_{1}=\frac{(2y_{5}y_{7}-y_{6}^{2})^{7}}{y_{7}^{12}}$ \\
& $[ V_{1},Y_{4}^{\prime}]=4Y_{4}^{\prime}, $ & $I_{2}=\frac{(3y_{7}^{2}(y_{4}^{\prime}-y_{4})-y_{6}^{3}+3y_{5}y_{6}y_{7})^{7}}{y_{7}^{18}}$\\
& $[ Y_{1},Y_{i}]  =Y_{i+1},\quad i=3,4,5,6 $ & \\
& $[ Y_{3},Y_{4}]=[Y_{3},Y_{4}^{\prime}]=Y_{7}.$ &\\\hline
$\frak{r}_{8}^{15}$ & $[ V_{1},Y_{i}]=iY_{i},\quad i=1,2,3,4,5,$ & $I_{1}=\frac{(2y_{3}^{\prime}y_{5}+y_{4}^{2}-2y_{4}y_{4}^{\prime})^{5}}{y_{5}^{8}}$ \\
& $[ V_{1},Y_{3}^{\prime}]=3Y_{3}^{\prime},\quad [V_{1},Y_{4}^{\prime}]  =4Y_{4}^{\prime},$ & $I_{2}=\frac{(y_{4}-y_{4}^{\prime})^{5}}{y_{5}^{4}}$\\
& $[ Y_{1},Y_{i}]  =Y_{i+1},\quad i=2,3,4,$ & \\
& $[ Y_{1},Y_{3}^{\prime}]=Y_{4}^{\prime},[Y_{1},Y_{4}^{\prime}]=[Y_{2},Y_{3}]=Y_{5}$ & \\\hline
$\frak{r}_{8}^{16}$ & $[ V_{1},Y_{i}]  =iY_{i},\quad i=1,2,3,4,5,$ & $I_{1}=\frac{(2y_{3}^{\prime}y_{5}-y_{4}^{2}+2y_{3}y_{5})^{5}}{y_{5}^{8}}$ \\
&$[ V_{1},Y_{1}^{\prime}]=Y_{1}^{\prime},\quad [V_{1},Y_{3}^{\prime}]  =3Y_{3}^{\prime},$ & $I_{2}=\frac{(y_{1}^{\prime}y_{5}-y_{2}y_{4}+y_{3}^{\prime}y_{3}-y_{1}y_{5})^{5}}{y_{5}^{6}}$ \\
& $[Y_{1},Y_{i}]  =Y_{i+1},\quad i=2,4,$ & \\
& $[ Y_{1},Y_{i}^{\prime}]=Y_{i+1},\quad i=1,3,$ & \\
&$[ Y_{1}^{\prime},Y_{i}]  =Y_{i+1},\quad i=3,4,$ & \\
& $[ Y_{1}^{\prime},Y_{2}]  =Y_{3}^{\prime}, $ & \\
& $[ Y_{2},Y_{3}]=-[Y_{2},Y_{3}^{\prime}]=Y_{5}$ & \\\hline
$\frak{r}_{8}^{17}$ &  $[ V_{1},Y_{i}]=iY_{i},\quad i=1,2,3,4,5,$ & $I_{1}=\frac{y_{5}^{\prime}}{y_{5}}$\\
& $[ V_{1},Y_{i}^{\prime}]  =iY_{i}^{\prime},\quad i=3,5,$ & $I_{2}=\frac{(2y_{3}y_{5}^{2}-2y_{3}^{\prime}y_{5}^{\prime}y_{5}-y_{4}^{2}(y_{5}^{\prime}-y_{5})+4y_{3}y_{5}y_{5}^{\prime})^{5}}{y_{5}^{8}(y_{5}^{\prime})^{5}}$ \\
&$[ Y_{1},Y_{i}]=Y_{i+1},\quad i=2,3,4,$ & \\
&$[ Y_{1},Y_{3}^{\prime}]=Y_{4},\quad [Y_{2},Y_{3}]=Y_{5}^{\prime},$ & \\
&$[ Y_{2},Y_{3}^{\prime}]  =Y_{5}+2Y_{5}^{\prime}$ & \\\hline
$\frak{r}_{8}^{18}$ & $[ V_{1},Y_{i}]=iY_{i},\quad i=1,2,3$ & $I_{1}=\frac{y_{3}^{\prime}}{y_{3}}$ \\
& $[V_{1},Y_{i}^{\prime}]=iY_{i}^{\prime},\quad i=1,2,3$ & $I_{2}=\frac{(-y_{1}^{"}(y_{3}^{\prime})^{2}+y_{3}^{\prime}(y_{2}^{2}+(y_{2}^{\prime})^{2})-2y_{1}(y_{3}^{2}-(y_{3}^{\prime})^{2})+2y_{3}(y_{1}^{\prime}y_{3}^{\prime}-y_{2}^{\prime}y_{2}))^{3}}{y_{3}(y_{3}^{\prime})^{6}}$\\
& $[V_{1},Y_{1}^{"}]=Y_{1}^{"}, \quad [Y_{1},Y_{1}^{\prime}]  =Y_{2},$ & \\
& $[ Y_{1},Y_{1}^{"}]  =Y_{2}^{^{\prime}},\quad [ Y_{1}^{"},Y_{2}^{\prime}]  =Y_{3}^{\prime},$ & \\
& $[ Y_{1},Y_{2}^{\prime}]=Y_{3}^{^{\prime}},\quad
[Y_{1}^{\prime},Y_{2}]=Y_{3}^{^{\prime}},$ & \\
& $[ Y_{1}^{\prime},Y_{2}^{\prime}]  =[Y_{1}^{"},Y_{2}^{{}}]  =Y_{3}$ & \\\hline
\bf
\end{tabular}
\end{table}

\begin{table}
\caption{\label{rigid82}Invariants of solvable rigid Lie algebras in dimension $8$ and rank $2$.}
\begin{tabular}{@{}lll}
\bf
Algebra&{\rm Brackets}&{\rm Invariants}\\\hline
\rm
$\frak{r}_{8}^{19}$ & $[V_{1},Y_{i}]=iY_{i},\ i=1,2,3,4,5,6,$ & $I_{1}=\frac{(y_{5}^{2}-2y_{4}y_{6})^{3}}{(3y_{6}^{2}y_{3}+y_{5}^{3}-3y_{4}y_{5}y_{6})^{2}}$\\
& $[V_{2},Y_{i}]=Y_{i},\quad i=2,3,4,5,6,$ & $I_{2}=\frac{(-6y_{2}y_{6}^{3}+6y_{5}y_{3}y_{6}^{2}+2y_{5}^{4}-8y_{4}y_{5}^{2}y_{6}+6y_{4}^{2}y_{6}^{2})}{(2y_{4}y_{6}-y_{5}^{2})^{2}}$\\
& $[Y_{1},Y_{i}]=Y_{i+1},\quad i=2,3,4,5.$ & \\\hline
$\frak{r}_{8}^{20}$ & $[V_{1},Y_{i}]=iY_{i},\quad i=1,2,3,4,5,6,$ & none\\
& $[V_{2},Y_{i}]=Y_{i},\quad i=3,4,5,6,$ &\\
& $[ Y_{1},Y_{i}]=Y_{i+1},\quad i=3,4,5,$ & \\
& $[ Y_{2},Y_{i}]  =Y_{i+2},\quad i=3,4$ & \\\hline
$\frak{r}_{8}^{21}$ &$[V_{1},Y_{i}]=iY_{i},\quad i=1,2,3,4,5,6,$ & none \\
& $[V_{2},Y_{i}]  =Y_{i},\quad i=4,5,6,$ & \\
& $[Y_{1},Y_{i}]  =Y_{i+1},\quad i=2,4,5,\quad [Y_{2},Y_{4}]  =Y_{6}$ & \\\hline
$\frak{r}_{8}^{22}$ &$[V_{1},Y_{i}]=iY_{i},\quad i=1,2,3,4,5,6,$ & none \\
& $[V_{2},Y_{i}]=Y_{i}, i=2,3,4,\quad [V_{2},Y_{i}]=2Y_{i}, i=5,6,$ & \\
& $[ Y_{1},Y_{i}]=Y_{i+1},\quad i=2,3,5,$ & \\
& $[ Y_{2},Y_{i}]=Y_{i+2},i=3,4$ & \\\hline
$\frak{r}_{8}^{23}$ &$[V_{1},Y_{i}]=iY_{i},\quad i=1,2,3,4,5,6,$ & none \\
& $[V_{2},Y_{6}]=Y_{6},\quad [Y_{2},Y_{3}]  =Y_{5}$ &\\
& $[ Y_{1},Y_{i}]  =Y_{i+1},\quad i=2,3,4,$ & \\\hline
$\frak{r}_{8}^{24}$ &$[V_{1},Y_{i}]=iY_{i},\quad i=1,2,3,4,5,7,$ & none \\
& $[V_{2},Y_{i}]=Y_{i},\quad i=3,4,5,\quad [V_{2},Y_{7}]  =2Y_{7}$ & \\
& $[Y_{1},Y_{i}]  =Y_{i+1},\quad i=3,4,$ & \\
& $[Y_{2},Y_{3}]=Y_{5},\quad [Y_{3},Y_{4}]  =Y_{7}$ & \\\hline
$\frak{r}_{8}^{25}$ &$[V_{1},Y_{i}]=iY_{i},\quad i=1,2,3,4,5,7,$ & none \\
& $[V_{2},Y_{i}]=Y_{i},\quad i=2,3,4,\quad [V_{2},Y_{5}]=2Y_{5},$ & \\
& $[V_{2},Y_{7}]=3Y_{7},\quad [Y_{1},Y_{i}]=Y_{i+1},\quad i=2,3,$ &\\
& $[Y_{2},Y_{i}]=Y_{i+2},\quad i=3,5$ & \\\hline
$\frak{r}_{8}^{26}$ &$[V_{1},Y_{i}]=iY_{i},\quad i=1,2,3,4,5,7,$ & none \\
& $[V_{2},Y_{i}]=Y_{i},\quad i=2,3,4,5,\quad [V_{2},Y_{7}]  =2Y_{7},$ & \\
& $[Y_{1},Y_{i}]  =Y_{i+1},\quad i=2,3,4,$ & \\
& $[Y_{2},Y_{5}]=-[Y_{3},Y_{4}]=Y_{7}$ & \\\hline
$\frak{r}_{8}^{27}$ &$[V_{1},Y_{i}]=iY_{i},\quad i=1,2,3,5,6,$& none \\
& $[V_{2},Y_{i}]=Y_{i},\quad i=2,3,\quad [V_{2},Y_{5}]=2Y_{5}, [V_{2},Y_{6}]=3Y_{6}$ & \\
& $[V_{2},Y_{7}]=3Y_{7}, [Y_{1},Y_{i}]  =Y_{i+1},\quad i=2,6,$& \\
& $[ Y_{2},Y_{i}]  =Y_{i+2},\quad i=3,5$ & \\\hline
$\frak{r}_{8}^{28}$ &$[V_{1},Y_{i}]=iY_{i},\quad i=1,2,3,4,5,7,$& none \\
& $[V_{2},Y_{2}]=Y_{2}, [V_{2},Y_{i}]=2Y_{i},\quad i=3,4,$ & \\
& $[V_{2},Y_{5}]=3Y_{5}, [V_{2},Y_{7}]  =4Y_{7}, [Y_{1},Y_{3}]=y_{4}$ & \\
& $[ Y_{3},Y_{4}]=-Y_{7}, [Y_{2},Y_{3}]=Y_{5}, [Y_{2},Y_{5}]=Y_{7} $ &\\
& $[Y_{2},Y_{i}]  =Y_{i+2},\quad i=3,5,$ &\\\hline
$\frak{r}_{8}^{29}$ &$[V_{1},Y_{i}]=iY_{i},\quad i=1,2,3,4,5,$& none \\
& $[ V_{1},Y_{3}^{\prime}]=3Y_{3}^{\prime}, [V_{2},Y_{3}^{\prime}]=Y_{3}^{\prime}, [V_{2},Y_{5}]=2Y_{5}$ &\\
& $[V_{2},Y_{i}]=Y_{i},\quad i=2,3,4,\quad [Y_{1},Y_{2}]=Y_{3}$ & \\
& $[Y_{1},Y_{3}]=Y_{4}, [Y_{1},Y_{3}^{\prime}]=Y_{4}, [Y_{2},Y_{3}]=Y_{5}$ & \\\hline
\bf
\end{tabular}
\end{table}

\begin{table}
\caption{\label{rigid83}Invariants of solvable rigid Lie algebras in dimension $8$ and rank $3$.}
\begin{tabular}{@{}lll}
\bf
Algebra&{\rm Brackets}&{\rm Invariants}\\\hline
\rm
$\frak{r}_{8}^{30}$& $[V_{1},Y_{i}]=iY_{i},\ i=1,2,3,4,5,$ & none\\
& $[V_{2},Y_{5}]=Y_{5},\quad [V_{3},Y_{i}]=Y_{i},\quad i=2,3,4,$ & \\
& $[ Y_{1},Y_{i}]=Y_{i+1},\quad i=2,3$ & \\\hline
$\frak{r}_{8}^{31}$ & $[V_{1},Y_{i}]=iY_{i},\quad i=1,2,3,4,5,$ & none\\
&$[V_{2},Y_{i}]=Y_{i},\quad i=2,3,$ & \\
& $[V_{3},Y_{i}]=Y_{i},\quad i=4,5,$ &\\
& $[Y_{1},Y_{i}]  =Y_{i+1},\quad i=2,4$ &\\\hline
$\frak{r}_{8}^{32}$ & $[V_{1},Y_{i}]=iY_{i},\quad i=1,2,3,4,5,$ & $I_{1}=\frac{v_{1}y_{5}+y_{1}y_{4}}{y_{5}}$\\
& $[V_{2},Y_{i}]=Y_{i},\quad i=2,4,5,$ & $I_{2}=\frac{v2y_{5}-v_{3}y_{5}+y_{2}y_{3}}{y_{5}}$\\
& $[V_{3},Y_{i}]  =Y_{i},\quad i=3,4,5,$ & \\
& $[ Y_{1},Y_{4}] =[Y_{2},Y_{3}]  =Y_{5}$ &\\\hline
\bf
\end{tabular}
\end{table}

The tables for rigid Lie algebras of rank at least two show a rather interesting fact, namely, that the algebra has nontrivial invariants if and only if it is isomorphic to either $\frak{r}_{8}^{19}$ or $\frak{r}_{8}^{32}$. If we analyze the nilradicals of these algebras, we see that the one corresponding to $\frak{r}_{8}^{19}$ is isomorphic to the filiform Lie algebra $L_{6}$ \cite{Ve}, while the second is isomorphic to the five dimensional Heisenberg Lie algebra $\frak{h}_{2}$. The interest of this is that these nilpotent Lie algebras are "extreme", in the sense that $L_{6}$ corresponds to the simplest algebra of maximal nipotence index, while the Heisenberg algebra is the most simple metabelian Lie algebra, which are the "most nilpotent". There is another important observation: while the rank of $\frak{h}_{p}$ increases with $p$, the rank of $L_{n}$ does not depend on the dimension, and equals $2$ \cite{Ve}.  

\begin{table}
\caption{\label{rigid84}Invariants of solvable rigid Lie algebras in dimension $8$ and rank $4$.}
\begin{tabular}{@{}lll}
\bf
Algebra&{\rm Brackets}&{\rm Invariants}\\\hline
\rm
$\frak{r}_{8}^{33}$& $[V_{1},Y_{i}]=iY_{i},\quad i=1,2,3,4,$ & none\\
& $[V_{2},Y_{2}]=Y_{2}, \quad [V_{3},Y_{3}]=Y_{3},\quad [V_{4},Y_{4}]  =Y_{4}$ & \\\hline
\bf
\end{tabular}
\end{table}

\section{Higher dimensions}
Although the root theory introduced in \cite{AG1} allows a complete classification of solvable rigid laws in dimension nine, this has not been established due to the great number of isomorphism classes, but only for particular cases \cite{AC,Car2}. In contrast to other solvable Lie algebras, for the rigid case we know that the torus, which is maximal, in fact determines the law on the nilradical. This information can be used to derive formulae of the invariants for arbitrary dimension, when the structure of the torus and its influence on the structure of the solutions in low dimension is known. We will illustrate this method with the following example: for $m\geq 2$ let $\frak{d}_{2m+1}=\frak{g}_{2m+1}\oplus\frak{t}$ be the solvable Lie algebra whose brackets over the basis $\left\{X_{0},..,X_{3},Y_{1},..,Y_{2m-3},V_{1},..V_{m-1}\right\}$ are 
\[
\begin{array}[c]{l}
[X_{0},X_{i}]=X_{i+1}, i=1,2;\quad [X_{1},Y_{2m-3}]=X_{3}\\
\lbrack Y_{2i-1},Y_{2i}]=X_{3}, i=1,..,m-2;\quad [V_{1},X_{i}]=X_{i},\quad i=0,2\\
\lbrack V_{1},X_{3}]=2X_{3};\quad [V_{1},Y_{2m-3}]=2Y_{2m-3}\\
\lbrack V_{1},Y_{2i}]=-2Y_{2i},\quad i=1,..,m-2\\
\lbrack V_{2},X_{i}]=X_{i}, i=1,2,3;\quad [V_{2},Y_{2i}]=Y_{2i},\quad 1\leq i\leq m-2\\
\lbrack V_{i+2},Y_{2i-1}]=Y_{2i-1}, i=1,..,m-2,\quad [V_{i+2},Y_{2i}]=-Y_{2i},\quad i=1,..,m-2 
\end{array}
\]
Using the root system associated to this algebra, it is easily seen that it is rigid. For $m=2$ it is isomorphic to $\frak{r}_{7}^{7}$, and expressed in the basis above, a fundamental set of invariants of this algebra is given by $\left\{\frac{-2x_{3}^{2}v_{2}+x_{3}^{2}v_{1}+x_{0}x_{1}x_{3}-2x_{1}x_{3}y_{1}+y_{1}x_{2}^{2}}{x_{3}^{2}}\right\}$. Now observe that, for $m\geq 3$, the nilradical of $\frak{d}_{2m+1}$ consists of a "germ subalgebra" generated by $\left\{X_{0},..,X_{3},Y_{1}\right\}$, to which three dimensional Heisenberg Lie algebras whose derived subalgebra is the center $\left(X_{3}\right)$ are "glued". We claim

\begin{proposition}
For any $m\geq 3$ the solvable Lie algebra $\frak{d}_{2m+1}$ has a fundamental set of invariants formed by the $(m-1)$ rational invariants given by
\[
\left\{\frac{-2x_{3}^{2}v_{2}+x_{3}^{2}v_{1}+x_{0}x_{1}x_{3}-2x_{1}x_{3}y_{2m-3}+y_{2m-3}x_{2}^{2}}{x_{3}^{2}}, \frac{x_{3}v_{k+2}+y_{2k-1}y_{2k}}{x_{3}}\right\}
\]
for $1\leq k\leq m-2$.
\end{proposition}

\begin{proof}
Over the coordinate system $\left\{x_{0},..,x_{3},y_{1},..,y_{2m-3},v_{1},..,v_{m}\right\}$ the system is the following:
\begin{eqnarray}
\widetilde{X}_{0}.F=\left(-x_{2}\partial_{x_{1}}-x_{3}\partial_{x_{2}}+x_{0}\partial_{v_{1}}\right).F=0 \\
\widetilde{X}_{1}.F=\left(x_{2}\partial_{x_{0}}-x_{3}\partial_{y_{2m-3}}+x_{1}\partial_{v_{2}}\right).F=0\\
\widetilde{X}_{2}.F=\left(x_{3}\partial_{x_{0}}+x_{2}\partial_{v_{1}}+x_{2}\partial_{v_{2}}\right).F=0\\
\widetilde{X}_{3}.F=\left(2x_{3}\partial_{x_{0}}+x_{3}\partial_{v_{2}}\right).F=0\\
\widetilde{Y}_{2i-1}.F=\left(-x_{3}\partial_{y_{2i}}+y_{2i-1}\partial_{v_{i+1}}\right).F=0,\quad 1\leq i\leq m-2\\
\widetilde{Y}_{2m}.F=\left(x_{3}\partial_{y_{2i-1}}+y_{2i}\left(2\partial_{v_{1}}+\partial_{v_{2}}-\partial_{v_{i+2}}\right)\right).F=0,\quad 1\leq i\leq m-2\\
\widetilde{Y}_{2m-3}.F=\left(x_{3}\partial_{x_{1}}+2y_{2m-3}\partial_{v_{1}}\right).F=0
\end{eqnarray}

Equations $(8.1),(8.2),(8.3),(8.4)$ and $(8.7)$ correspond to the system of partial differential equations associated to the seven dimensional Lie algebra $\frak{r}_{7}^{7}$, with solution $\frac{-2x_{3}^{2}v_{2}+x_{3}^{2}v_{1}+x_{0}x_{1}x_{3}-2x_{1}x_{3}y_{2m-3}+y_{2m-3}x_{2}^{2}}{x_{3}^{2}}$. Now this function is easily seen to be a solution of the system above. Equations $(8.5)$ and $(8.6)$ have, for fixed $i$, the solution  $\frac{x_{3}v_{k+2}+y_{2k-1}y_{2k}}{x_{3}}$, which is also a solution of the whole system. Observe that these solution have all the same structure, since they correspond to the infinitesimal generators $\widetilde{Y}_{2i-1}, \widetilde{Y}_{2i}$ of the Heisenberg subalgebras generated by them, jointly with the element $\widetilde{X}_{3}$. The rank of the commutator table is easily seen to be $m-1$.
\end{proof}

\section{Conclusions}

We have computed the invariants of complex solvable rigid Lie algebras up to dimension eight. This is indeed the maximal dimension for which these algebras have been classified, although certain partial classifications exist. Most of the general results obtained for arbitrary dimension also remain valid in the real case. However, there is a reason for studying complex rigid Lie algebras rather than real ones. It is known that that any complex simple Lie algebra is in fact defined over the rational $\mathbb{Q}$, and from the tables above it also follows that any solvable complex rigid Lie algebra is also rational. Since we have seen that the eigenvalues for the action of a maximal torus of derivations can be chosen in $\mathbb{Z}$, Carles \cite{Car} conjectured that any rigid Lie algebra is rational. This would justify in an elegant manner why simple Lie algebras are rational, since they are rigid by cohomological reasons. Unfortunately this conjecture is false, as pointed out by Goze and Ancochea in \cite{AG3}. They gave examples of rigid Lie algebras being not only non-rational, but even non-real. That these algebras remained unknown relies on the fact that this pathology begins to appear from dimension eleven on.\newline We have pointed out that, in the rigid case, useful information on the existence or not of non-trivial invariants can be deduced from the action of the torus on its nilradical, without paying closer attention to the structure of the latter. This works since this structure is perfectly determined by the torus. The example in section 8. also shows that we can determine a fundamental set of invariants starting from solutions in low dimensional cases. This will be of particular importance when we consider semidirect products of nilpotent Lie algebras and non-maximal tori. Although in this case the algebras are not rigid, probably we can deduce information about the invariants by analyzing the action of these tori. This approach should provide at least sufficiency criteria for the non-existence of non-trivial invariants, since the obtention of characterizations for these algebras probably do not exist. Such sufficiency conditions will be of extreme importance for the non-decomposable case.\newline The study of invariants of solvable rigid Lie algebras also seems to answer another important question: it is well known that for semisimple Lie algebras the cardinal of a fundamental set of invariants coincides with the rank of the algebra, i.e., the dimension of a Cartan subalgebra \cite{Ra}. It is therefore natural to ask wheter there exists a link between the cardinal of a fundamental set and the rank for the solvable case. If we consider the $n+3$ dimensional solvable Lie algebra $L_{n}\oplus\frak{t}$ given by 
\begin{eqnarray*}
\lbrack X_{1},X_{i}]=X_{i+1},\quad i=2,..,n\\
\lbrack V_{1},X_{i}]=i.X_{i}, \quad i=1,..,n+1\\
\lbrack V_{2},X_{i}]=X_{i}, \quad i=2,..,n+1
\end{eqnarray*}
we see that the rank is always two, and a Cartan subalgebra is generated by $V_{1}$ and $V_{2}$. Now the number of functionally independent invariants for this algebra is $\mathcal{N}=n-3$, and therefore depends on the dimension of the nilradical, and not on the dimension of a Cartan subalgebra. We may expect that for more complicated examples, such as non-decomposable solvable Lie algebras,  such a formula will generally be not inferable. \newline Finally, we have seen by examples that a solvable Lie algebra having a fundamental set of invariants formed by rational functions needs not to be rigid. However, there is a wide class of algebras having this property, namely those satisfying $dim \frak{g}=dim Der\left(\frak{g}\right)$, where $Der\left(\frak{g}\right)$ denotes the Lie algebra of derivations. Specifically, this is of interest for complete Lie algebras \cite{AC2,Z}, i.e., centerless Lie algebras whose derivations are all inner. This class is of great importance for structure theory, and its invariants will surely provide interesting information about their representations.

\end{document}